\documentclass[letterpaper, reqno]{amsart}
\usepackage{amssymb, amsmath, amsthm}
\usepackage[numbers]{natbib}
\usepackage[all]{xy}
\usepackage{hyperref}
\newcommand{\V}{\mathcal{V}}
\newcommand{\hh}{\mathcal{H}}
\newcommand{\In}{\mathcal{I}}

\newcommand{\R}{\mathbb{R}}

\newcommand{\Ii}{\frak{I}}
\newcommand{\X}{\frak{X}}
\newcommand{\inv}{^{\text{-}1}}

\newcommand{\vf}{\frak{X}}

\newcommand{\del}[1]{\frac{\partial}{\partial #1}}
\newcommand{\dev}[1]{\frac{d}{d #1}}
\newcommand{\Del}[1]{\frac{D}{\partial #1}}

\swapnumbers
\numberwithin{equation}{section} 

\newtheorem{theorem}{Theorem}[section]
\newtheorem{lem}[theorem]{Lemma}

\newtheorem{proposition}[theorem]{Proposition}
\newtheorem{cor}[theorem]{Corollary}
\newtheorem{co}[theorem]{Conjecture}
\newtheorem*{ta}{Theorem}

\theoremstyle{definition}
\newtheorem{df}[theorem]{Definition}

\theoremstyle{remark}
\newtheorem{rem}[theorem]{Remark}
\newtheorem{que}{Question}

\newcommand{\bd}{\begin{df}}
\newcommand{\ed}{\end{df}}

\newcommand{\bq}{\begin{que}}
\newcommand{\eq}{\end{que}}

\newcommand{\bl}{\begin{lem}}
\newcommand{\el}{\end{lem}}

\newcommand{\br}{\begin{rem}}
\newcommand{\er}{\end{rem}}

\newcommand{\bt}{\begin{theorem}}
\newcommand{\et}{\end{theorem}}

\newcommand{\bc}{\begin{cor}}
\newcommand{\ec}{\end{cor}}

\newcommand{\bco}{\begin{co}}
\newcommand{\eco}{\end{co}}

\newcommand{\bp}{\begin{proposition}}
\newcommand{\ep}{\end{proposition}}

\newcommand{\brm}{\begin{rem}}
\newcommand{\erm}{\end{rem}}

\newcommand{\be}{\begin{equation}}
\newcommand{\ee}{\end{equation}}

\begin{document}
\title{Submanifolds and the Sasaki metric}
\author{Pedro Sol\'orzano}
\begin{abstract}This is the content of a talk given by the author at the 2009 Lehigh University Geometry/Topology Conference. 

Using the definition of connection given by Dieudonn\'e, the Sasaki metric on the tangent bundle to a Riemannian manifold is expressed in a natural way.

Also, the following property is established.   The induced metric on the tangent bundle of an isometrically embedded submanifold is the Sasaki metric if and only if the submanifold is totally geodesic.
\end{abstract}

\maketitle
\section{Introduction}
The second order tangent bundle to a Riemannian manifold comes with a splitting that facilitates the construction of Riemannian metrics on its total space. These metrics (often called {\em natural}) have arisen in different settings (cf. \cite{MR946027}, \cite{MR2106375}) and usually inherit many properties from the geometry of their base manifold. 

The Sasaki metric (introduced in \cite{MR0112152}) is a standard pick, although \citet{MR946027} proved it is quite rigid. This rigidity can be seen as a particular case of a more general interplay between the metric on $TM$ and on $M$, as noted by \citet{MR2106375}. 

\citet{MR946027} also introduced another natural metric which they called the {\em Cheeger-Gromoll metric} (the first explicit computation of a metric suggested by \citet{MR0309010} for the tangent bundle of the sphere); they studied these metrics when restricted to the unit tangent bundle as well. Both the tangent bundle and the unit tangent bundle have very rich structures, the former as an almost complex manifold and the latter as contact manifold. The geometry of the unit tangent bundle is coupled with that of the manifold with many a conjecture (cf. \cite{MR1798731}). All these directions have also proven prolific, as presented in the survey by \citet{MR1798731}.

The Sasaki metric occurs in a natural way as a Whitney sum metric on a vector bundle by virtue of Proposition \ref{natmet}. Furthermore, the heritability of the Sasaki metrics under smooth mappings is studied. In particular, the following result is verified (presented here in \ref{tgs}).
\begin{ta} The induced metric on total space of the tangent bundle of an isometrically imbedded submanifold (from the Sasaki metric on the total space of the tangent bundle of the ambient manifold) coincides with its Sasaki metric if and only if the embedding is totally geodesic.
\end{ta}
This is a consequence of Theorem \ref{pullbk}, which relates the two Sasaki metrics with the second fundamental form of the embedding.

\section{Vector Bundles}
\subsection*{Background}
Let us recall some classical constructions and introduce some notation. Every diagram is assumed to be commutative unless otherwise stated.

\bd\label{whitsum}Let $(E,\pi_1)$ and $(F,\pi_2)$ be vector bundles over a manifold $M$. Define their {\em Whitney sum} as a vector bundle over $M$ with total space denoted by $E\oplus F$ and projection map $\pi_1\circ pr_1=\pi_2\circ pr_2$ fitting into the universal diagram for the pullback.  Namely,

\be
\xymatrix{
E\oplus F 
\ar[r]^{pr_2}
\ar[d]_{pr_1}
& F\ar[d]_{\pi_2}
 \\
E
\ar[r]^{\pi_1}
& M
}
\ee
\ed

\br
The projections from the Whitney sum to its factors are also vector bundles regarding, e.g., the pullback $\pi_1^*=E\oplus$  as a functor from the category of bundles over $M$ to the category of bundles over $E$. 
\er 

\bp\label{2struc} 
Given a vector bundle, $(E,\pi)$, $\pi:E\rightarrow M$, there are two vector bundle structures with total space $TE$, namely the standard projection, $(TE,\pi_E)$,
$$\pi_E:TE\rightarrow E,$$
and the  {\em secondary structure}, $(TE,\pi_*)$, 
$$\pi_*:TE\rightarrow TM.$$
Let $TE\oplus_2 TE$ denote the Whitney sum using the secondary structure; that is $(X,Y)\in TE\oplus_2TE$ satisfy $\pi_*X=\pi_*Y$:

\be
\xymatrix{
TE\oplus_2 TE 
\ar[r]^{ \phantom{aaa}pr_2}
\ar[d]_{pr_1}
& TE\ar[d]_{\pi_*}
 \\
TE
\ar[r]^{\pi_*}
& TM
}
\ee
\ep

In the case when $(E,\pi)=(TM,\pi_M)$, the base space is $TM$ for both structures and one has the following result.

\bp[\citet{MR1724021}]\label{invol} There exists a bundle involution $\In:TTM\rightarrow TTM$ from one structure to another. Namely, one has the following commutative diagram

\be
\xymatrix{
TTM 
\ar[rr]^{\In}
\ar[dr]_{\pi_{TM}}
&& TTM
\ar[dl]^{{\pi_{M}}_*}
 \\
&TM
}
\ee
\ep
%\br 
%In local co\"ordinates the definition of $\In$ is clear.
%\er

\bd
Following the notation by \citet*{MR0365399}, given a (normed) vector space $V$, there is a canonical isomorphism between $V\times V$ and $TV$, given by
\be\label{curli}
\Ii_v(w)f=\Ii(v,w)f=\dev{t}\bigg|_{t=0}f(v+tw).
\ee
That is, $\Ii_vw$ is the directional derivative at $v$ in the direction $w$. 
\ed

\bp Given any vector bundle $(E,\pi)$, (\ref{curli}) yields a bundle isomorphism between $\oplus^2E:=E\oplus E$ and the vertical distribution $\V=\ker\pi_*\subseteq TE$, in a natural way; that is, there is a natural transformation (also denoted by $\Ii$) from the functor $\oplus^2$ to the functor $T$.
\ep
\begin{proof}
The naturality: Let $(E,\pi_1)$ and $(F,\pi_2)$ be vector bundles over $M$, and let $\varphi:E\rightarrow F$ be a morphism between them. Then,
\be
\xymatrix{
\oplus^2 E 
\ar[r]^{\oplus^2\varphi}
\ar[d]_{\Ii}
& \oplus^2F
\ar[d]_{\Ii}
 \\
TE
\ar[r]^{\varphi_*}
& TF
}  
\phantom{aaaaaa}
\xymatrix{
(e,\tilde e)
\ar@{|->}[r]
\ar@{|->}[d]
& (\varphi e,\varphi\tilde e)
\ar@{|->}[d]
 \\
[e+t\tilde e]_e
\ar@{|->}[r]
&\varphi_*([e+t\tilde e]_e)=[\varphi_* e+t\varphi_* \tilde e]_{\varphi_* e}
}  
\ee
with $[\alpha(t)]_{\alpha(0)}=\dot\alpha(0)$, where $\alpha$ is a curve. The fact that it maps into the vertical distribution follows from 
\be
{\pi_1}_*[e+t\tilde e]_e=[\pi_1(e+t\tilde e)]_{\pi e}=[\pi_1 e]_{\pi e}=0,
\ee 
since by assumption $\pi_1 e=\pi_1\tilde e=\pi_1(e+t\tilde e)$. Surjectivity can also be verified.
\end{proof}

\bc\label{vstar}
Let $f:M\rightarrow N$ be a smooth map between smooth manifolds. Then
\be
f_{**}\circ\Ii=\Ii\circ(\oplus^2 f_*).
\ee
\ec

\bd[\citet{MR0362066}] A connection on a vector bundle $(E,\pi)$ is a bundle morphism $C:E\oplus TM\rightarrow TE$ with respect to both bundle structures on $TE$:
\be
\xymatrix{
E\oplus TM 
\ar[rr]^{C}
\ar[dr]_{pr_1}
&& TE
\ar[dl]^{\pi_E}
 \\
&E
}
\phantom{aaaaaa}
\xymatrix{
E\oplus TM 
\ar[rr]^{C}
\ar[dr]_{pr_2}
&& TE
\ar[dl]^{\pi_*}
 \\
&TM
}
\ee
\ed

\br\label{connint} 
One should read $C(e,u)$ as the ``horizontal lift of $u$ at $e$", since given a connection $C$ one can define the horizontal space as $\hh_e=C(\{e\}\times T_{\pi(e)}M)$ as well as a projection onto the vertical space $\V$, also denoted by $\V$.
\er

\br\label{covdev}
The standard (Kozul) covariant derivative definition of a connection is equivalent and is recovered by the following equation. Let $Y:M\rightarrow E$ be a section of the bundle $(E,\pi)$, let $x\in TM$; then
\be\label{covdeveq}
\Ii(Y,\nabla_xY)=Y_*x-C(Y,x).
\ee 
\er
\subsection*{Bundle Isomorphisms}

Let us now give a few different presentations of $TE$ and $TTM$, that will provide the setting for \ref{natmet}.
%The following proposition might prove useful in the sequel.

\bp\label{3map}
Given a connection $C$ on $(E,\pi)$ there exists a bundle isomorphism $\Xi=\Xi_C:E\oplus TM\oplus E\rightarrow TE$ as bundles over $E$. 
\ep
\begin{proof} Define $\Xi$ by 
\be\label{3m}
\Xi(e,u,f)=C(e,u)+\Ii(e,f).
\ee
This map is a bundle map in view of the following diagram.
\[
\xymatrix{
E\oplus TM\oplus E 
\ar[rr]^{\phantom{aaaa}\Xi}
\ar[dr]_{pr_1}
&& TE
\ar[dl]^{\pi_E}
 \\
&E
}
\phantom{aaaaaaa}
\xymatrix{
(e,u,f)
\ar@{|->}[rr]
\ar@{|->}[dr]
&& C(e,u)+\Ii(e,f)
\ar@{|->}[dl]
 \\
&e
}
\]

 In order to prove that it is an isomorphism, an inverse can be produced:
 \be
 \Xi\inv(X)=(\pi_EX,\pi_*X,\Ii_{\pi_EX}\inv\V X)
 \ee
with $\V$ as in \ref{connint}.
\end{proof}

\bc\label{3map} 
$TTM$ is bundle isomorphic to $\oplus^3TM$.
\ec

\bp\label{2vert}
Given a connection $C$ on $TM$, $TTM$ is bundle isomorphic to $\oplus^2\V$ over $TM$. 
\ep
\begin{proof}
By \ref{3map}, one has only to check that $\oplus^2\V$ is isomorphic to $\oplus^3TM$.  Since $\V\cong\oplus^2 TM$, 
\be
\oplus^2(\oplus^2 TM)=\{(X,Y)|pr_1 X=pr_1Y\}=\{((u, v),(u,w))\}.
\ee
Thus the isomorphism is given by $(u,v,u,w)\mapsto(u,v,w)$. 
\end{proof}

\subsection*{Metric structures}
We now review the concept of metric structures on vector bundles and the work of \citet{MR1724021}.

\bd
A metric on a vector bundle $(E,\pi)$ is a function $g:\oplus^2E\rightarrow\R$ satisfying the usual conditions. Given $(E,\pi)$ and a vector bundle with metric $(F, \tilde\pi,h)$ there is a natural metric on $\pi^*F=E\oplus F$ as a bundle over $E$ given by the pullback metric
\be
\pi^*h=h\circ(\oplus^2pr_2).
\ee
\ed

\br\label{bunmet}
Given two bundles with metrics $(E, \pi_1,g)$, $(F, \pi_2,h)$ over $M$, there is a natural metric on their Whitney sum as bundles over $M$:
\be
g\oplus h=g\circ(\oplus^2pr_1)+h\circ(\oplus^2pr_2).
\ee
\qed
\er

\citet{MR1724021} present a characterization of the Levi-Civita connection in this language. For completeness, let us recall some of their results.

\bd[\citet{MR1724021}] Let $g:\oplus^2E\rightarrow\R$ be a metric on a vector bundle $(E,\pi)$. There is an associated map $\widetilde{Tg}:TE\oplus_2 TE\rightarrow\R$ given by
\be
\widetilde{Tg}(X,Y)=\dev{t}\bigg|_{t=0}g(u(t),v(t)),
\ee
where $u$ and $v$ are curves in $E$ with $\pi\circ u=\pi\circ v$ such that $\dot u(0)=X$ and $\dot v(0)=Y$ (Recall $\pi_*X=\pi_*Y$ from \ref{2struc}).
\ed

Applying these definitions to $E=TM$, one has:
\bp[\citet{MR1724021}]\label{LCcheck} Given a connection $C$ on $(TM,\pi_M)$,

it is metric iff
\be \widetilde{Tg}(C(u,h),C(v,h))=0 \text{ }\forall h,u,v\in TM;
\ee 

and it is torsion free iff
\be C(u,v)=\In C(v,u)\text{ }\forall u,v\in TM
\ee
\ep

\section{Sasaki Metric on $TM$}

Let us review the Sasaki metric and prove a proposition demonstrating that the Sasaki metric may be viewed as a Whitney sum metric.

Several metrics on $TM$ have been studied. An important class is that of metrics which render $\pi_M$ a Riemannian submersion while preserving the natural splitting on $TTM$ given in \ref{tsplit}; these are called $g$-natural metrics and have been studied profusely (c.f \citet{MR2106375}).

\bd\label{tsplit}
Let $(M,g)$ be a Riemannian manifold. The metric induces a natural splitting of $(TTM,\pi_{TM})$ as a bundle over $TM$ into vertical and horizontal subbundles, as noted in \ref{connint}:

\be
\hh=C(\oplus^2TM),
\ee

\be
\V=\Ii(\oplus^2TM).
\ee
\ed

\br
By virtue of \ref{covdev}, $\hh$ can also be realized as $\{\dot X| \Del{t}X=0\}$, where $X$ is a vector field along some curve in $M$, $\Del{t}$ is the covariant derivative on that curve, and $\dot X= X_*\del{t}$. That is
\be
\hh=\{\dot X| X\text{ parallel along }\pi\circ X\}.
\ee

\er

\bd 
For a given vector field $X\in\vf(M)$ there are two associated vector fields $X^v,X^h\in\vf(TM)$ given by
\be
X^v(u)=\Ii(u,X),
\ee
\be
X^h(u)=C(u,X).
\ee
These are called the {\em vertical} and {\em horizontal lifts} of $X$.  It is standard to check that they are smooth.
\ed

With this splitting in mind, one can define several metrics on $TM$ that turn $\pi$ into a Riemannian submersion. 

\bd
The {\em Sasaki metric} on $TM$, denoted by $g_S$, is defined by requiring vertical and horizontal vectors be orthogonal to each other and by 

\be\label{sasmet}
g_S(X^v,Y^v)=g(X,Y)\circ\pi=g_S(X^h,Y^h).
\ee
\ed

\bp \label{natmet}
Let $\Psi:TTM\rightarrow\oplus^2\V$ be the bundle isomorphism guaranteed by \ref{2vert}. Let 
\be
h=(pr_2\circ\Ii\inv)^*g\oplus(pr_2\circ\Ii\inv)^*g
\ee
be the metric on $\oplus^2\V$ (cf. \ref{bunmet}). Then $\Psi$ is a bundle isometry.
\ep
\begin{proof}
At a point $u\in TM$, (\ref{sasmet}) can be written as follows

\be\label{altsas}
\begin{array}{c}
g_S\circ (\oplus^2C_u)=g\circ (\oplus^2pr_2)=g_S\circ (\oplus^2\Ii_u)\\ \\
g_S\circ(\Ii_u\oplus C_u)=0,
\end{array}
\ee
%by recalling that 
%\be
%\xy
%(-15,8)*+{\oplus^2TM}="t0"; 
%(-15,-8)*+{TM}="b0";
%(15,8)*+{TM}="t1"; 
%(15,-8)*+{M}="b1";
%{\ar^{pr_2} "t0"; "t1"};
%{\ar^{ \pi_M} "b0"; "b1"};
%{\ar^{ pr_1} "t0"; "b0"};
%{\ar^{\pi_M} "t1"; "b1"};
%\endxy 
%\Longrightarrow
%\xy
%(-15,8)*+{\oplus^2(\oplus^2TM)}="t0"; (-15,-8)*+{\oplus^2TM}="b0";
%(15,8)*+{\oplus^2 TTM}="t1"; (15,-8)*+{\R}="b1";
%{\ar@<.5ex>^{\oplus^2\Ii_u} "t0"; "t1"};
%{\ar@<-.5ex>_{\oplus^2C_u} "t0"; "t1"};
%{\ar^{ g} "b0"; "b1"};
%{\ar^{ \oplus^2pr_2} "t0"; "b0"};
%{\ar^{g_S} "t1"; "b1"};
%\endxy 
%\ee
with $pr_2$ as in \ref{whitsum}. The claim now follows by looking at \ref{bunmet} (noticing that in the current setting there is a systematic abuse of notation as to the meaning of $pr_2$).
\end{proof}

\section{Submanifolds}
It is natural to ask when is the Sasaki metric preserved under mappings between Riemannian manifolds. 

Consider the case when $\iota:(M,g)\rightarrow(\bar M,\bar g)$ is an isometric immersion, i.e. $g=\iota^*\bar g$. Take a point $x\in TM$ and a vector field $y\in\X(M)$. Let $Y\in\X(\bar M)$ be a vector field $\iota$-related to $y$, i.e. 
\be\label{rlted}
Y\circ\iota=\iota_*\circ y.
\ee

\bd
Define the {\em second fundamental form} $B:\oplus^2TM\rightarrow T\bar M$ by
\be
B(x,y)=\bar\nabla_{\iota_*x}Y-\iota_*\nabla_xy.
\ee
%It is standard to check that it is well defined, bilinear and symmetric.
\ed

\bl\label{vsas}
The relation between the vertical parts at a point $u\in TM$ is encoded by the following equation.
\be 
(\iota_*)^*\bar g_S\circ (\oplus^2\Ii_u)=g_S\circ (\oplus^2\Ii_u)
\ee
\el
\begin{proof} 
For the second equality \ref{vstar} is used:
\be
\text{LHS}:=\bar g_S\circ\iota_{**}\circ (\oplus^2\Ii_u)=(\bar g_S\circ (\oplus^2\Ii_{\iota_* u}))\circ(\oplus^2\iota_*)=\bar g\circ (\oplus^2\iota_*)=g=:\text{RHS.}
\ee

\end{proof} 

\bl\label{cstar}
Let $C$ and $\bar C$ be their corresponding Levi-Civita connections. Then,
\be
\iota_{**}C=\bar C\circ (\oplus^2\iota_*)+\Ii\circ ((\iota_*\circ pr_1)\oplus B).
\ee
\el
\begin{proof} 
In view of \ref{vstar}, one has that
\be
\iota_{**}\circ\Ii=\Ii\circ(\oplus^2 \iota_*).
\ee
which combined with (\ref{covdeveq}) produces
\begin{align}
\iota_{**}\Ii(y,\nabla_xy) &= \Ii(\iota_*y,\bar\nabla_{\iota_*x}Y-B(y,x))\nonumber\\
&=Y_*\iota_*x-\bar C(Y\circ\iota,\iota_*x)-\Ii(\iota_*y,B(y,x))
\end{align}
and
\be
\iota_{**}\Ii(y,\nabla_xy)=\iota_{**}y_*x-\iota_{**}C(y,x).
\ee
Solving for $\iota_{**}C(y,x)$ and using (\ref{rlted}) and its derivative version yields
 \be
\iota_{**}C(y,x)=\bar C(\iota_*y,\iota_*x)+\Ii(\iota_*y,B(y,x)).
\ee

\end{proof}

\bt \label{pullbk}Let  $\iota:(M,g)\rightarrow(\bar M,\bar g)$ is an isometric immersion at a point $u\in TM$. Then
\be\label{pllbkeqn}
(\iota_*)^*\bar{g}_S=g_S+\bar g\circ(\oplus^2B_u).
\ee
\et
\begin{proof} 
By virtue of \ref{vsas}, \ref{cstar}, and (\ref{altsas}) a direct computation yields the equation
\be
(\iota_*)^*\bar{g}_S(X,Y)=g_S(X,Y)+\bar g(B(u,X),B(u,Y)),
\ee
where $X,Y\in T_uTM$.
\end{proof}

\bc\label{tgs}
The induced metric coincides with the Sasaki metric iff the embedding is totally geodesic.
\ec

\br The previous statement clearly hold in greater generality. More details about this fact will follow.
\er

\bibliographystyle{plainnat}
\bibliography{ref}
\end{document}